\documentclass[11pt]{article}
\usepackage[a4paper,margin=28mm]{geometry}
\usepackage[T1]{fontenc}
\usepackage{lmodern,microtype,amsmath,amssymb,amsthm,mathtools,tikz}
\usepackage[colorlinks=true,linkcolor=black,citecolor=black,urlcolor=blue]{hyperref}
\usetikzlibrary{arrows.meta}
\hypersetup{pdftitle={Rigidity and nonrigidity for interval translation maps},pdfauthor={Pascal Hubert and Alexandra Skripchenko}}
\numberwithin{equation}{section}
\newtheorem{theorem}{Theorem}[section]
\newtheorem{proposition}[theorem]{Proposition}
\newtheorem{lemma}[theorem]{Lemma}

\theoremstyle{definition}
\theoremstyle{remark}\newtheorem{remark}[theorem]{Remark}
\newtheorem*{namedremark}{Remark}
\newcommand{\R}{\mathbb R}
\newcommand{\Q}{\mathbb Q}
\newcommand{\Z}{\mathbb Z}
\newcommand{\N}{\mathbb N}
\newcommand{\one}{\boldsymbol{1}}
\newcommand{\ab}{\operatorname{ab}}

\newcommand{\diam}{\operatorname{diam}}
\newcommand{\calA}{\mathcal A}
\newcommand{\calR}{\mathcal R}
\newcommand{\calC}{\mathcal C}
\newcommand{\calE}{\mathcal E}
\newcommand{\calP}{\mathcal P}
\newcommand{\eps}{\varepsilon}
\title{Rigidity and nonrigidity for interval translation maps}
\author{Pascal Hubert \and Alexandra Skripchenko}
\date{}
\begin{document}
\maketitle
\begin{abstract}
Almost every infinite-type interval translation map in the
Bruin--Troubetzkoy family, and in Bruin's family with any fixed number
of branches, is rigid for its invariant probability, where typicality
is understood with respect to the intrinsic induction parameter measure.
We prove this by exhibiting explicit repeated-word blocks in the
induction and showing that their cylinders recur.
In the opposite direction, we give a sufficient condition for
nonrigidity of a primitive aperiodic substitution: the existence of
a rational periodic vector of its incidence matrix in the return-word
space whose coordinate sum is nonzero. The proof gives a uniform positive
lower bound for symbol discrepancies at all sufficiently large times.
For self-similar Bruin maps with primitive period substitution, the
same condition also implies weak mixing. We apply the nonrigidity
condition to an explicit self-similar infinite-type Bruin map on eight
intervals, specifying its algebraic parameters and verifying its
return words.
\end{abstract}

\section{Introduction and main results}
Interval translation maps (ITMs) are maps of an interval that act by
translations on the elements of a finite interval partition. They were
introduced by Boshernitzan and Kornfeld \cite{BK} as a generalization
of interval exchange transformations (IETs), in which the images of
the continuity intervals are allowed to overlap. Their renormalization,
attractors and ergodic properties have been actively studied; see, for
example, \cite{AAFHS,AHS,BR,DSS,Kryzhevich,CS}. Much of the existing work
concerns specific families of ITMs, notably the Bruin--Troubetzkoy family
\cite{BT}, Bruin's family \cite{Bruin} and double rotations \cite{BC}.
In this paper we continue the study of the ergodic properties of the
Bruin--Troubetzkoy and Bruin families by examining rigidity.

Our first theorem establishes unique ergodicity and rigidity for
almost every nondegenerate infinite-type map in the Bruin--Troubetzkoy
family and in Bruin's family with any fixed number of branches, with
respect to the intrinsic induction parameter probability. Our second
theorem gives a general sufficient condition for nonrigidity of
primitive aperiodic substitution systems. Our third theorem applies
this condition to construct an explicit self-similar Bruin map on
eight intervals whose attracting dynamics is minimal, uniquely ergodic
and nonrigid.

There are two different measures in these statements. The intrinsic
parameter probability describes typical maps in a family. For each
individual map, rigidity concerns its invariant probability on the
attractor. In particular, almost all infinite-type maps does not
mean almost all parameters for ambient Lebesgue measure.

A probability-preserving map $(Y,\mu,T)$ is \emph{rigid} if there are
integers $q_j\to\infty$ such that
\begin{equation}\label{eq:rigidity}
 \|f\circ T^{q_j}-f\|_{L^2(\mu)}\longrightarrow0
 \qquad\text{for every }f\in L^2(\mu).
\end{equation}
Equivalently, $\mu(T^{-q_j}A\triangle A)\to0$ for every measurable
set $A$. The sequence is the same for all observables; the convergence
is strong operator convergence, not convergence in operator norm.
Rigidity for general IETs was first studied by Veech
\cite{VeechRigidity}, who proved that Lebesgue-almost every IET is rigid.

Fix $d\ge2$, put $N=d+1$, and let
\[
 U_d=\{1>\alpha_1>\cdots>\alpha_d>0\}.
\]
For $\alpha\in U_d$, set $\alpha_0=1$ and define Bruin's map by
\begin{equation}\label{eq:itm}
 T_\alpha(x)=
 \begin{cases}
 x+\alpha_i,&x\in[1-\alpha_{i-1},1-\alpha_i),\quad 1\le i\le d,\\
 x+\alpha_d-1,&x\in[1-\alpha_d,1).
 \end{cases}
\end{equation}
Let $I_1,\ldots,I_N$ be the half-open intervals in \eqref{eq:itm},
in order from left to right. We write $\calP=\{I_1,\ldots,I_N\}$
for the continuity partition and $\calP(x)=i$ for the label of
the interval containing $x$.
The case $d=2$ is the Bruin--Troubetzkoy family \cite{BT}. An interval translation
map $T$ is of \emph{finite type} if $T^{n+1}[0,1)=T^n[0,1)$ for some
$n$, and of \emph{infinite type} otherwise.

The Gauss induction \cite[Section~2, equation~(1)]{Bruin} is
\begin{equation}\label{eq:gauss}
 G_d(\alpha)=
 \left(\frac{\alpha_2}{\alpha_1},\ldots,\frac{\alpha_d}{\alpha_1},
 k+\frac{\alpha_d-1}{\alpha_1}\right),
 \qquad k=\left\lfloor\frac1{\alpha_1}\right\rfloor .
\end{equation}
Write $\calE_d=\bigcap_{n\ge0}G_d^{-n}(U_d)$ for the nondegenerate
survivor set. Thus all $N$ branches are preserved by every induction
step. For larger $N$, parameters which first lose branches and
subsequently have infinite-type dynamics in a smaller family are
not included in $\calE_d$.

Let $\mathbb P_d$ be the intrinsic induction parameter probability on
$\calE_d$: for $d=2$ it is the law constructed by
Artigiani--Hubert--Skripchenko \cite{AHS}, and for general $d$ it is
the law of Skripchenko--Chernyi \cite{CS}. We specify its coding
convention and the recurrence input in Sections~\ref{sec:parameter}
and~\ref{bruin:sec:parameter}.
These measures are supported on sets of zero ambient Lebesgue measure.

\begin{theorem}\label{thm:typical}
For every $d\ge2$, there is a Borel subset of $\calE_d$ of full
$\mathbb P_d$-measure such that each corresponding map $T_\alpha$
has a unique invariant probability $\mu_\alpha$ and is rigid with
respect to $\mu_\alpha$.
\end{theorem}

The proof of Theorem~\ref{thm:typical} uses a standard argument based on repeated return words and
follows, for example, the approach of Ferenczi--Hubert \cite{FH} to
rigidity for interval exchange transformations. The substitution and
tower framework is classical \cite{DHS}.

Our main nonrigidity result, Theorem~\ref{thm:criterion}, is a general
sufficient condition for primitive aperiodic substitutions and thus
applies beyond the setting of ITMs. Its application to ITMs here is
restricted to self-similar maps, for which a return of the induction
yields a substitution. The proof compares every candidate
time with a suitable substitution scale. A signed combination of actual
return lengths produces a fixed nonzero displacement, which cannot be
an approximate period of an aperiodic linearly recurrent language.

We next state a sufficient condition for nonrigidity. For a substitution $\sigma$
on a finite alphabet $\calA$, use the column-count convention
\[
 M_{ij}=|\sigma(j)|_i,\qquad
 \ab(\sigma(w))=M\ab(w),\qquad h_n=\one M^n.
\]
Here $\ab(w)$ is a column vector, $\one=(1,\ldots,1)$ is a row,
and $h_n(i)=|\sigma^n(i)|$. A \emph{return word to a letter $a$}
starts with $a$, has no other occurrence of $a$, and is followed by
$a$ in an admissible word. The \emph{return space} is
\[
 \calR_\sigma=\operatorname{span}_{\R}
 \{\ab(r): r\text{ is a return word to some letter}\}.
\]
This definition requires actual admissible words; closed paths in the
graph of admissible pairs need not themselves be admissible words.

\begin{theorem}\label{thm:criterion}
Let $\sigma$ be a primitive aperiodic substitution with incidence
matrix $M$ and invariant probability $\nu$. Suppose that there are
$p\ge1$ and $u\in\Q^{\calA}$ such that
\begin{equation}\label{eq:criterion}
 M^p u=u,\qquad \one u\ne0,\qquad u\in\calR_\sigma.
\end{equation}
Then there exist $\eps>0$ and $q_0\ge1$ such that
\begin{equation}\label{eq:gap}
 \nu\{x:x_0\ne x_q\}\ge\eps\qquad(q\ge q_0).
\end{equation}
In particular, $(X_\sigma,\nu,S)$ is not rigid.
\end{theorem}

For self-similar Bruin maps with primitive period substitution, this
nonrigidity condition also implies weak mixing; see
Lemma~\ref{lem:bruin-nonrigidity-weak-mixing}.

\begin{namedremark}
Following Mercat \cite[Section~2.3]{Mercat}, a real letter-return
coboundary is a row $b$ such that $b\ab(r)=0$ for every return word
$r$. These rows form $\calC_\sigma=\calR_\sigma^\perp$.
The construction in \cite[Lemma~2.4 and its proof]{Mercat} gives
a function $g$, depending only on the initial letter, with
$b_{x_0}=g(x)-g(Sx)$. Thus $u\in\calR_\sigma$ is equivalent to
$bu=0$ for every $b\in\calC_\sigma$. The coboundaries are linear
forms on letter-count columns; the rational column $u$ in
\eqref{eq:criterion} belongs to their common kernel. In particular,
the return-space condition is automatic when $\calC_\sigma=0$.
\end{namedremark}

\begin{theorem}\label{thm:eight}
Let $\varrho$ be the unique root in $(0,1/3)$ of
\begin{equation}\label{eq:Q}
 Q(t)=t^7-2t^6-2t^3+2t^2-4t+1.
\end{equation}
Define
\begin{align}
 a_7&=\frac{1-3\varrho}{1-\varrho},&
 a_3&=\frac{a_7}{\varrho},&
 a_4&=\varrho,\notag\\
 a_1&=1-a_7+\varrho^2,&
 a_2&=1-a_7+\varrho^2a_1,&
 a_5&=\varrho a_1,\qquad a_6=\varrho a_2.
 \label{eq:parameters}
\end{align}
Then $1>a_1>\cdots>a_7>0$, and the eight-interval map
$T_a$ defined by \eqref{eq:itm} is of infinite type.
Its attracting dynamics is minimal and uniquely ergodic and is not
rigid. Its Gauss digits are $(1,1,1,3)^\infty$, and its return
substitution is
\begin{equation}\label{eq:zeta}
 \begin{split}
 \zeta:\quad&1\mapsto5,\quad2\mapsto6,\quad3\mapsto7,\quad4\mapsto81,\\
 &5\mapsto82,\quad6\mapsto83,\quad7\mapsto8444,\quad8\mapsto844.
 \end{split}
\end{equation}
For the invariant probability $\mu$ and continuity partition $\calP$,
there exists a constant $\eps_8>0$ such that
\begin{equation}\label{eq:explicitgap}
 \mu\{x:\calP(x)\ne\calP(T_a^q x)\}
 \ge \eps_8\qquad(q\ge6).
\end{equation}
\end{theorem}

Artigiani--Avila--Ferenczi--Hubert--Skripchenko
\cite[Theorems~2 and 3]{AAFHS} prove that, for every fixed number
of branches, almost every Bruin interval translation map of infinite
type is weakly mixing with respect to its invariant probability.
Here typicality is taken with respect to the intrinsic parameter
law $\mathbb P_d$ used in Theorem~\ref{thm:typical}.
Intersecting the two full-measure sets therefore shows that
$\mathbb P_d$-almost every map in $\calE_d$ is both rigid and weakly mixing.

\begin{remark}[An obstruction to our nonrigidity construction in the BT family]
\label{rem:obstruction}
The obstruction below concerns the construction of a nonrigid
self-similar example by Theorem~\ref{thm:criterion} within the
three-interval Bruin--Troubetzkoy family. It follows from the spectral
analysis of self-similar interval translation maps: see
Boshernitzan--Kornfeld \cite{BK} for the original infinite-type
construction and Artigiani--Avila--Ferenczi--Hubert--Skripchenko
\cite[Proposition~13]{AAFHS} for the general periodic statement.
Indeed, let $P$ be a primitive period matrix of a self-similar BT map.
The latter result gives a strictly positive second Lyapunov exponent.
After telescoping to a positive period matrix, periodicity and
$|\det P|=1$ imply that two eigenvalues of $P$ have modulus greater
than one and the third has modulus less than one. Consequently
\[
 \det(P^s-I)\ne0
 \quad\text{and}\quad
 \det(P^s+I)\ne0
 \qquad(s\ge1).
\]
In particular, $P$ has no nonzero rational periodic vector, so the
relation $P^p u=u$ required by Theorem~\ref{thm:criterion} cannot hold.
This is an obstruction to the present construction of a nonrigid
BT example; it does not exclude another nonrigidity construction.
\end{remark}

\section{Repeated return words and typical rigidity}\label{sec:rigidity}

\subsection{Structure of the proof}

We first consider the Bruin--Troubetzkoy family and then extend the
technique to Bruin's family with any fixed number of branches.
The proof has two parts: a deterministic estimate for a prescribed
block of induction digits, and an almost-sure recurrence argument
with respect to the induction parameter measure.

The deterministic argument follows the repeated-return-word method
of Ferenczi--Hubert \cite{FH}. At a suitable induction depth, every
tower name contains a long power of the same earlier return word
$U$. Translation by $q=|U|$ preserves the continuity label except
on a small proportion of each tower. Weighting these proportions
by the tower masses gives a discrepancy estimate for every
invariant probability. Since the partitions into finite continuity itineraries
generate the Borel $\sigma$-algebra on $[0,1)$
(Lemma~\ref{lem:generating}), the estimate yields strong $L^2$ rigidity. The same repetitions
make the cylinder frequencies independent of the invariant probability.

For BT maps the block $(2,K,2)$ gives
\[
 \mu\{x:\calP(T^q x)\ne\calP(x)\}\le\frac9{K+7},
 \qquad q=|\theta_s(2)|,
\]
where $s$ is the position before the block and $\theta_s$ is the
substitution product defined below. The first and last digits $2$
provide, respectively, the required height comparison and a long
repeated word in every tower. Section~\ref{sec:parameter} shows that
these blocks occur with arbitrarily large $K$ for almost every
parameter under the Artigiani--Hubert--Skripchenko measure.
For Bruin maps, only the length of the strings of digits $2$, the
height calculation and the corresponding induction cylinder change.

Throughout the deterministic proof, $\mu$ is an invariant probability
of one fixed interval map, $\nu$ its symbolic image and
$\widehat\nu$ its two-sided extension. The parameter
measure $\mathbb P_d$ is used only for recurrence of induction blocks.

\subsection{The three-interval Bruin--Troubetzkoy family}
\label{bt:sec:proof}

For technical simplicity, we first treat the three-interval
Bruin--Troubetzkoy family, namely the case \(d=2\) of
\eqref{eq:itm}. The extension to Bruin's family is given in
Section~\ref{bruin:sec:extension}. Write \(I=[0,1)\), let
\(T=T_{\alpha,\beta}\), and label its continuity intervals by
\[
 I_1=[0,1-\alpha),\qquad
 I_2=[1-\alpha,1-\beta),\qquad
 I_3=[1-\beta,1).
\]
Thus \(\calP(x)=i\) for \(x\in I_i\).  Put
\(\Omega=\bigcap_{n\ge0}T^nI\).  An invariant probability \(\mu\) on \(I\)
has \(\mu(T^nI)=1\), since \(T^{-n}(T^nI)=I\), and consequently
\(\mu(\Omega)=1\).  All probability estimates below concern such invariant
probabilities, not ambient interval length.

For the BT Gauss digits \((k_n)\), the return substitutions and their
column-count matrices are
\begin{equation}\label{bt:eq:sub}
 \chi_k(1)=2,\qquad
 \chi_k(2)=31^k,\qquad
 \chi_k(3)=31^{k-1},\qquad
 M_k=
 \begin{pmatrix}
 0&k&k-1\\
 1&0&0\\
 0&1&1
 \end{pmatrix}.
\end{equation}
Set \(\theta_0=\mathrm{id}\) and
\(\theta_n=\chi_{k_1}\circ\cdots\circ\chi_{k_n}\).  The height row is
\[
 h_n=\bigl(|\theta_n(1)|,|\theta_n(2)|,|\theta_n(3)|\bigr),
 \qquad h_{n+1}=h_nM_{k_{n+1}}.
\]
The word \(\theta_n(i)\) records an actual return excursion at induction
depth \(n\), in the original continuity labels; its length counts
iterates of the original map.

\begin{theorem}\label{bt:thm:det}
Let \(T\) be an infinite-type BT map whose directive sequence contains
\[
 (k_{s_j+1},k_{s_j+2},k_{s_j+3})=(2,K_j,2),
 \qquad K_j\longrightarrow\infty,
\]
at increasing positions.  Then \(T\) has a unique invariant probability
\(\mu\) and is rigid with respect to \(\mu\).
After passing to disjoint blocks with \(K_j\ge2\), rigidity times are
\begin{equation}\label{bt:eq:times}
 q_j=|\theta_{s_j}\chi_2(1)|=|\theta_{s_j}(2)|.
\end{equation}
For \(m\ge1\), write
\(\calP_m(x)=(\calP(x),\calP(Tx),\ldots,\calP(T^{m-1}x))\)
for the length-\(m\) continuity itinerary of \(x\). Then
\begin{equation}\label{bt:eq:finite-name-bound}
 \mu\{x:\calP_m(x)\ne\calP_m(T^{q_j}x)\}
 \le \min\left\{1,\frac{9m}{K_j+7}\right\}.
\end{equation}
\end{theorem}

We first construct the symbolic towers used to count discrepancies.
We then identify the common repeated word and estimate the discarded
levels.  The final steps explain why those estimates give rigidity
for all measurable observables and why the same repetitions imply
unique ergodicity.  The almost-everywhere occurrence of the digit
blocks is a separate parameter argument.

\subsubsection{Continuity coding and invariant probabilities}

Let \(\sigma\) be the one-sided shift on
\(\{1,2,3\}^{\mathbb N_0}\).  Under the block condition used below,
the prefix limit
\[
 \rho=\lim_{n\to\infty}\theta_n(3)
\]
exists; this and its tail versions are proved in
Lemma~\ref{bt:lem:growth}.  Define
\[
 X=\overline{\{\sigma^r\rho:r\ge0\}}.
\]
For a word \(w=w_0\cdots w_{m-1}\), write
\([w]=\{x\in X:x_0\cdots x_{m-1}=w\}\).
These are orbit-word cylinders for one fixed map, as opposed to
the parameter cylinders used later.

The following statement collects classical coding and word-frequency
arguments from symbolic dynamics, specialized to the BT induction.

\begin{proposition}\label{bt:prop:coding}\label{bt:prop:frequencies}
Suppose that the directive sequence of an infinite-type BT map
contains infinitely many blocks $(2,K,2)$ with $K\ge2$.
Its closed attractor $\overline\Omega$ is a Cantor set, the map
has no periodic points, and an invariant probability on $I$ exists.
The forward continuity coding
\[
 \iota(x)_r=\calP(T^rx),\qquad r\ge0,
\]
is Borel, satisfies $\iota\circ T=\sigma\circ\iota$, and takes
values in $X$ outside a countable endpoint exceptional set in
the attractor. Every invariant probability $\mu$ of $T$ is
nonatomic, and $\iota_*\mu$ is a shift-invariant probability on $X$.

At any occurrence of $(2,K,2)$ with $K\ge2$ after a prefix
$\theta_s$, put $U=\theta_s(2)$ and $a=|U|$.
For a word $w$ of length $m\ge1$, set
\[
 p_U(w)=\frac1a
 \#\{0\le r<a:(U^\infty)_{[r,r+m)}=w\}.
\]
Every invariant probability $\nu$ on $X$ satisfies
\begin{equation}\label{bt:eq:frequencies}
 \bigl|\nu([w])-p_U(w)\bigr|
 \le\frac{8a+m-1}{(K-1)a}.
\end{equation}
Consequently, the block condition of Theorem~\ref{bt:thm:det}
implies unique ergodicity of $X$.
\end{proposition}

\begin{proof}
The attractor and coding assertions follow from
\cite[Propositions~3--4, Corollary~5, Theorem~7 and the proof of
Lemma~17]{BT}. Existence of an invariant probability follows from
\cite[Corollary~4.15]{BR}, whose hypotheses hold under the recurring
block assumption.
The frequency estimate is the normalized return-word frequency
argument of \cite[Proposition~3.22 and its proof]{VeechIET}
and \cite[Lemma~17 and its proof]{BT}, applied to the tower names
\eqref{bt:eq:Z} and the length bounds \eqref{bt:eq:balance}
established below.
\end{proof}

\subsubsection{Growth and symbolic concatenations}

At a fixed induction depth, concatenating two return words means
following one excursion by the other.  The following growth statement
ensures both that the relevant infinite words exist and that the
candidate return times tend to infinity.

\begin{lemma}\label{bt:lem:growth}
If the directive sequence contains infinitely many disjoint blocks
\((2,K,2)\), with \(K\ge2\), then
\[
 \min_{1\le i\le3}|\theta_n(i)|\longrightarrow\infty.
\]
The limit \(\rho=\lim_n\theta_n(3)\) exists.  At every fixed depth
\(r\), it is a concatenation of
\(\theta_r(1)\), \(\theta_r(2)\) and \(\theta_r(3)\).
Every finite factor of every point of \(X\) occurs in \(\rho\).
\end{lemma}

\begin{proof}
This is the standard prefix-limit construction; see
\cite[Theorem~7 and its proof]{BT} and \cite[Proposition~5]{Bruin}.
Every $M_k$ has nonzero columns, so the minimum height is
nondecreasing. Moreover,
\begin{equation}\label{bt:eq:positive-block}
 \one M_2M_KM_2=(K+2,K+7,K+4)\ge4\one
 \qquad(K\ge2).
\end{equation}
Each selected block therefore multiplies the minimum height by
at least four, proving growth. Since every $\chi_k(3)$ begins
with $3$, the words $\theta_n(3)$ are nested prefixes.
The same argument applies to each tail of the directive sequence;
writing $\rho^{(r)}$ for its limit gives
$\rho=\theta_r(\rho^{(r)})$ and hence the required concatenation.
The assertion about factors follows from the definition of $X$
as the shift orbit closure of $\rho$.
\end{proof}

\subsubsection{Return-word towers}

We use the standard return-word tower construction, as in
Ferenczi \cite[Proposition~2.24 and the following paragraph]{FerencziInduction}.
For the BT substitutions, recognition of word boundaries is immediate
from their form. Let $\widehat X$ be the two-sided subshift with the
language of $\rho$, let $S$ be its shift, and let $\widehat\nu$ be the
shift-invariant two-sided extension of an invariant probability $\nu$ on $X$.
A tower of base $B$ and height $H$ consists of the disjoint levels
$B,SB,\ldots,S^{H-1}B$; its name records their continuity labels.

\begin{lemma}\label{bt:lem:towers}
Under the hypothesis of Lemma~\ref{bt:lem:growth}, every
\(y\in\widehat X\) has a unique decomposition into
\(\theta_n(1),\theta_n(2),\theta_n(3)\) at every depth \(n\).
Let \(B_n(i)\) be the set where coordinate zero is the beginning
of a depth-\(n\) word of type \(i\).  Then
\begin{equation}\label{bt:eq:KR-partition}
 \widehat X
 =\bigsqcup_{i=1}^3\bigsqcup_{r=0}^{h_n(i)-1}S^rB_n(i)
\end{equation}
is a clopen Kakutani--Rokhlin partition.  The \(i\)-th tower has
name \(\theta_n(i)\).  The decompositions at successive depths
are compatible: the next towers stack copies of the preceding
towers in the order specified by \(\chi_{k_{n+1}}\).
\end{lemma}

\begin{proof}
At a step with digit $k$, the non-$1$ symbols mark exactly the
beginnings of the images $2$, $31^k$ and $31^{k-1}$.
A $2$ identifies type $1$, while a $3$ followed by a maximal run
of $k$ or $k-1$ symbols $1$ identifies type $2$ or $3$.
This includes $k=1$. Every two-sided word in the tail language
admits this parsing: the runs of $1$ are bounded, and every finite
block between successive markers occurs in
$\rho^{(r)}=\chi_{k_{r+1}}(\rho^{(r+1)})$.
Including the next marker also shows that the parsed type sequence
belongs to the next tail language.

Iteration gives a unique type and offset at every depth, hence
\eqref{bt:eq:KR-partition}. For fixed depth, parsing uses finitely
many coordinates, so the bases and levels are clopen. The next
boundary after $B_n(i)$ occurs at time $h_n(i)$, with name
$\theta_n(i)$. Finally,
$\theta_{n+1}=\theta_n\chi_{k_{n+1}}$ gives the stacking rule.
\end{proof}

Put $\beta_n(i)=\widehat\nu(B_n(i))$. Invariance gives
\begin{equation}\label{eq:towermass}
 \widehat\nu(S^rB_n(i))=\beta_n(i),\qquad
 \sum_{i=1}^3h_n(i)\beta_n(i)=1.
\end{equation}
Consequently, discarding at most $\varepsilon h_n(i)$ levels from
each tower discards total probability at most $\varepsilon$,
irrespective of the base masses.

\subsubsection{The height cone and the common repeated word}

Fix an occurrence of \((2,K,2)\) after the prefix
\(\theta=\theta_s\), with \(K\ge2\).

\begin{lemma}\label{bt:lem:height-cone}
Every prefix height row \(h_n=(x,y,z)\) satisfies
\begin{equation}\label{bt:eq:height-cone}
 x\le y,\qquad z\le y,\qquad y\le x+z.
\end{equation}
\end{lemma}

\begin{proof}
The inequalities hold for \(h_0=(1,1,1)\).
The update is
\[
 (x,y,z)M_k=(y,kx+z,(k-1)x+z).
\]
For the new row, the first inequality follows from
\(y\le x+z\le kx+z\), the second is immediate, and the third
is equivalent to \(x\le y\).  Hence the cone is invariant
for every \(k\ge1\).
\end{proof}

After the first digit \(2\), put
\[
 U=\theta\chi_2(1),\quad
 V=\theta\chi_2(2),\quad
 W=\theta\chi_2(3),\qquad
 a=|U|,\quad b=|V|,\quad c=|W|.
\]
These are the names and heights of the depth-\((s+1)\) towers.
Since \(\chi_2(1)=2\), \(\chi_2(2)=311\), and \(\chi_2(3)=31\),
if \(h_s=(x,y,z)\), then
\begin{equation}\label{bt:eq:height-formulas}
 a=y,\qquad b=2x+z,\qquad c=x+z.
\end{equation}
The height cone gives
\begin{equation}\label{bt:eq:balance}
 c\le2a,\qquad b\le3a,\qquad c+2b\le8a.
\end{equation}
Thus one preliminary digit \(2\) controls the lengths needed
in the argument, uniformly in the preceding directive sequence.

The next two digits satisfy
\[
 \chi_K\chi_2(1)=31^K,\qquad
 \chi_K\chi_2(2)=31^{K-1}22,\qquad
 \chi_K\chi_2(3)=31^{K-1}2.
\]
Consequently the depth-\((s+3)\) tower names are exactly
\begin{equation}\label{bt:eq:Z}
 Z_1=WU^K,\qquad
 Z_2=WU^{K-1}V^2,\qquad
 Z_3=WU^{K-1}V.
\end{equation}
All three contain a long power of the same word \(U\).
The final digit \(2\) is needed for this conclusion:
after the digit \(K\) alone, the type-\(1\) name would be \(V\),
because \(\chi_K(1)=2\), and would have no designated repeated
\(U\)-run.

The candidate time is \(q=a=|U|\).
At offset \(t\), with \(0\le t<a\), in a copy of \(U\) followed
by another copy, shifting by \(a\) reaches the same offset in
the next copy.  The current continuity label is therefore
unchanged.  The two points occupy different levels of the larger
tower; what is used here is equality of their labels.
This is the repeated-return-word mechanism in the rigidity
argument of \cite{FH}.
For successive disjoint blocks it produces the times
\eqref{bt:eq:times}; Lemma~\ref{bt:lem:growth} implies that
these times tend to infinity.

\subsubsection{Discarded levels and the discrepancy estimate}

In each tower in \eqref{bt:eq:Z}, retain every copy of \(U\)
except the last.  Discard the prefix \(W\), the final copy of \(U\),
and the suffix made of \(V\)'s, if present.
Every retained level has its \(a\)-step comparison inside the same
tower and with the same continuity label.
The discarded and retained counts are, respectively,
\[
 \begin{aligned}
 d_1&=c+a,   &g_1&=(K-1)a,\\
 d_2&=c+a+2b,&g_2&=(K-2)a,\\
 d_3&=c+a+b, &g_3&=(K-2)a.
 \end{aligned}
\]
Here \(|Z_i|=d_i+g_i\).  The estimates
\eqref{bt:eq:balance} give
\[
 d_1\le3a,\qquad d_2\le9a,\qquad d_3\le6a,\qquad
 g_i\ge(K-2)a.
\]
Some discarded levels may also give matching labels; the argument
only needs an upper bound on the exceptional set.

\begin{figure}[htbp]
\centering
\begin{tikzpicture}[x=1cm,y=.60cm,>=Stealth,font=\small]
\filldraw[fill=gray!15,draw=black!65] (0,0) rectangle (2.6,1);
\node at (1.3,.5) {$W$};
\foreach \j in {1,2,3,4}{
  \filldraw[fill=green!14,draw=black!65] (0,\j) rectangle (2.6,\j+1);
  \node at (1.3,\j+.5) {$U$};
}
\filldraw[fill=orange!25,draw=black!65] (0,5) rectangle (2.6,6);
\node at (1.3,5.5) {last $U$};
\foreach \j in {6,7}{
  \filldraw[fill=gray!15,draw=black!65] (0,\j) rectangle (2.6,\j+1);
  \node at (1.3,\j+.5) {$V$};
}
\draw[blue!70!black,thick] (0,1.3)--(2.6,1.3);
\draw[blue!70!black,thick] (0,2.3)--(2.6,2.3);
\draw[->,blue!70!black,thick] (-.45,1.3)--(-.45,2.3)
  node[midway,left] {$S^a$};
\node[anchor=west,text width=7.2cm] at (3.15,6.8)
  {The $V$-suffix and final $U$ are discarded.};
\node[anchor=west,text width=7.2cm] at (3.15,3.2)
  {In each retained copy, $a=|U|$ iterates lead to the same offset
   in the next copy of $U$.};
\node[anchor=west,text width=7.2cm] at (3.15,.5)
  {The initial $W$ is discarded.};
\node[anchor=north] at (1.3,-.25) {base of the $Z_2$-tower};
\end{tikzpicture}
\caption{The name $Z_2=WU^{K-1}V^2$ for $K=6$, read from bottom
to top. Each rectangle represents a whole word block, containing
$|W|$, $|U|$, or $|V|$ individual levels. The drawing is schematic,
not to scale. Green blocks give certified matches under $S^a$.}
\label{fig:tower}
\end{figure}
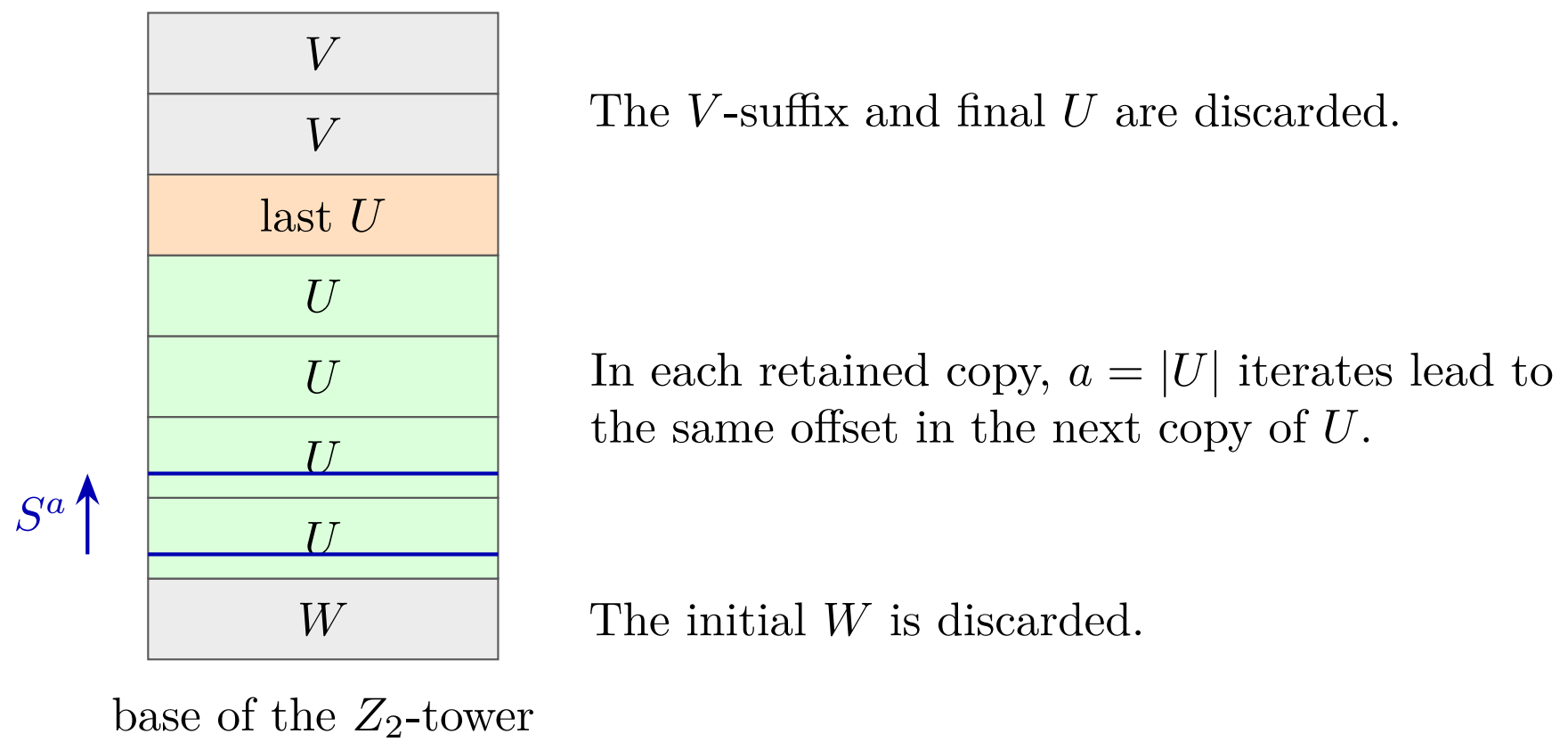

Figure~\ref{fig:tower} illustrates the retained copies in the second tower.

\begin{proposition}\label{bt:prop:mismatch}
At an occurrence of $(2,K,2)$ after $\theta_s$, with $K\ge2$, consider
the depth-$(s+3)$ tower partition with names $Z_i$ in \eqref{bt:eq:Z}.
Put $a=|\theta_s(2)|$.
Every invariant probability \(\nu\) on \(X\) satisfies
\begin{equation}\label{bt:eq:mismatch}
 \nu\{x:x_0\ne x_a\}\le\frac9{K+7}.
\end{equation}
\end{proposition}

\begin{proof}
The estimates $d_i\le9a$ and $g_i\ge(K-2)a$ give
$(K-2)d_i\le9g_i$. Since $|Z_i|=d_i+g_i$, this is equivalent to
\begin{equation}\label{bt:eq:bad-proportion}
 \frac{d_i}{|Z_i|}\le\frac9{K+7}.
\end{equation}
Let $B_i$ be the base of the $Z_i$-tower and put
$\beta_i=\widehat\nu(B_i)$, where $\widehat\nu$ is the two-sided
extension of $\nu$.
The discrepancy event is contained in the union \(E\) of
the discarded levels.  Using \eqref{eq:towermass},
\begin{equation}\label{bt:eq:weighted-bad}
 \widehat\nu(E)=\sum_{i=1}^3d_i\beta_i
 \le\frac9{K+7}\sum_{i=1}^3|Z_i|\beta_i
 =\frac9{K+7}.
\end{equation}
Projection to nonnegative coordinates proves \eqref{bt:eq:mismatch}.
\end{proof}

By Proposition~\ref{bt:prop:coding}, the proposition applies to
\(\nu=\iota_*\mu\) for every invariant probability of the interval map.
It follows that
\begin{equation}\label{bt:eq:interval-mismatch}
 \mu\{x:\calP(T^ax)\ne\calP(x)\}\le\frac9{K+7}.
\end{equation}

\subsubsection{From finite itineraries to rigidity}

We recall the two elementary facts that pass from a one-symbol
discrepancy estimate to rigidity. They apply to every finite-branch
interval translation map with a generating continuity partition.

\begin{lemma}\label{lem:generating}
Let \(R:I\to I\) be an interval translation map with finitely many
half-open continuity intervals, with continuity partition \(\mathcal Q\).
Suppose that
\(\overline{\bigcap_{n\ge0}R^nI}\) has empty interior.
Then the maximum diameter of an atom of
\[
 \mathcal Q_m=\bigvee_{r=0}^{m-1}R^{-r}\mathcal Q
\]
tends to zero.  In particular, forward continuity coding is
injective and generates the Borel sigma algebra on \(I\).
\end{lemma}

\begin{proof}
Set $C_n=\overline{R^nI}$. These compact sets are nested and
$C_n\setminus R^nI$ is finite. Thus $\bigcap_nC_n$ has empty
interior: otherwise an interval, after removal of countably many
endpoints, would lie in $\bigcap_nR^nI$, contrary to the hypothesis
on its closure. By compactness, the largest component diameter
$\eta_n$ of $C_n$ tends to zero; otherwise a sequence of components
would have a nondegenerate interval limit contained in every $C_n$.

Each atom $A$ of $\mathcal Q_m$ is an interval on which $R^{m-1}$
is a translation. Hence $\diam A\le\eta_{m-1}$.
Shrinking mesh gives injectivity of the coding and generation of
the Borel sigma algebra.
\end{proof}

\begin{lemma}\label{bt:lem:rigidity}
Let \(R\) preserve a Borel probability \(\mu\), and let
\(\mathcal Q\) be a finite partition whose forward itinerary
partitions generate modulo \(\mu\).
Suppose that \(q_j\to\infty\) and
\[
 \delta_j=\mu\{x:\mathcal Q(R^{q_j}x)\ne\mathcal Q(x)\}
 \longrightarrow0.
\]
Then
\(\|f\circ R^{q_j}-f\|_{L^2(\mu)}\to0\)
for every \(f\in L^2(\mu)\).
\end{lemma}

\begin{proof}
For a bounded function $g$ constant on the atoms of
$\mathcal Q_m=\bigvee_{r=0}^{m-1}R^{-r}\mathcal Q$, invariance gives
\[
 \|g\circ R^{q_j}-g\|_2^2
 \le4\|g\|_\infty^2m\delta_j\longrightarrow0,
\]
since disagreement of the first $m$ symbols is contained in the
union of the first $m$ inverse images of the one-symbol
disagreement set. Such functions are dense in $L^2(\mu)$, and
$f\mapsto f\circ R^{q_j}$ is an isometry. Approximation proves
the claim.
\end{proof}

\subsubsection{Completion of the deterministic BT proof}

\begin{proof}[Proof of Theorem~\ref{bt:thm:det}]
Choose disjoint occurrences with $K_j\to\infty$ and define $q_j$
by \eqref{bt:eq:times}. Lemma~\ref{bt:lem:growth} gives
$q_j\to\infty$. An invariant probability exists by
Proposition~\ref{bt:prop:coding}, and
Proposition~\ref{bt:prop:mismatch} gives, for every such probability,
\[
 \mu\{x:\calP(T^{q_j}x)\ne\calP(x)\}
 \le\frac9{K_j+7}\longrightarrow0.
\]
The Cantor attractor allows us to apply
Lemma~\ref{lem:generating}; Lemma~\ref{bt:lem:rigidity} then proves
rigidity and \eqref{bt:eq:finite-name-bound}.
The unique ergodicity of $X$ in
Proposition~\ref{bt:prop:coding} gives the same symbolic image for
all invariant interval probabilities. Since the continuity partition
generates, these probabilities coincide.
\end{proof}

\subsection{Recurrence of the required BT blocks}\label{sec:parameter}

The measure introduced by Artigiani--Hubert--Skripchenko
is ergodic and gives positive mass to every admissible cylinder
\cite[Section~3.1 and the proof of Lemma~19]{AHS}.
With the marked moves of \cite[Section~2.2 and Proposition~13]{AHS},
a complete Gauss digit $k$ is $A^{k-1}C_A$ or $B^{k-1}C_B$,
according to the initial state.

\begin{proposition}\label{bt:prop:parameter-recurrence}
For $\mathbb P_2$-almost every parameter, the Gauss sequence
contains $(2,2,2,K,2)$ infinitely often for each integer $K\ge2$.
\end{proposition}
\begin{proof}
For each $K\ge2$, the admissible closed path
\begin{equation}\label{bt:eq:synchronized-path}
 \mathsf W_K
 =C_A\,(BC_B)\,(AC_A)\,(BC_B)\,
       (A^{K-1}C_A)\,(BC_B)
\end{equation}
has a positive-measure cylinder. Its initial switch completes the
preceding digit, and the remaining runs give $(2,2,2,K,2)$.
Ergodicity and a countable intersection over $K$ prove the assertion.
\end{proof}

Choosing disjoint occurrences with $K_j\to\infty$ and taking their
last three digits, Theorem~\ref{bt:thm:det} proves the BT case of
Theorem~\ref{thm:typical}.

\subsection{Extension to Bruin's family}
\label{bruin:sec:extension}

The tower-mass calculation, the frequency argument and the passage
to $L^2$ rigidity apply to any finite alphabet. We give only the
additional combinatorial calculations and the induction cylinder
needed for Bruin's family with $N=d+1$ branches.

\subsubsection{Substitutions and the required block}

The return substitutions in \cite[Section~4, equation~(5)]{Bruin},
with the alphabet relabelled as $\{1,\ldots,N\}$, are
\begin{equation}\label{eq:chi}
 \chi_k(i)=i+1\quad(1\le i<d),\qquad
 \chi_k(d)=N1^k,\qquad \chi_k(N)=N1^{k-1},
 \qquad N=d+1.
\end{equation}
By \cite[Proposition~5]{Bruin}, the map has an attracting Cantor set
whose dynamics is identified with the corresponding substitution system.
The correspondence of invariant probabilities is described in
\cite[Lemma~6 and its proof]{Bruin}.
The endpoint exception is treated as in
Proposition~\ref{bt:prop:coding}: minimality of the Cantor attracting
dynamics excludes periodic points, so invariant probabilities are
nonatomic and give no mass to the countable exceptional set.
We retain the notation $\theta_s$, $X_d$ and $\widehat X_d$ for the
substitution products and their one-sided and two-sided languages,
with prefix limit $\rho_d=\lim_n\theta_n(N)$.

The block replacing $(2,K,2)$ is
\begin{equation}\label{eq:block}
 (\underbrace{2,\ldots,2}_{d-1},K,
   \underbrace{2,\ldots,2}_{d-1}).
\end{equation}
Set $A_d=d(d+7)/2$.

\begin{proposition}\label{prop:block}
Let $\alpha\in\calE_d$ have disjoint occurrences of
\eqref{eq:block} after prefixes of lengths $s_j$, with
$K_j\ge2$ and $K_j\to\infty$. Then $T_\alpha$ has a unique invariant
probability $\mu$ and is rigid for this probability. The times
$q_j=|\theta_{s_j}(d)|$ tend to infinity. For every fixed $m\ge1$,
\begin{equation}\label{eq:mainupper}
 \mu\{x:\calP_m(x)\ne\calP_m(T_\alpha^{q_j}x)\}
 \le\min\left\{1,\frac{mA_d}{K_j+A_d-2}\right\},
\end{equation}
where $\calP_m$ is the first-$m$ continuity name. The case $d=2$
gives the BT estimate.
\end{proposition}

\subsubsection{The height comparison and repeated-word calculation}

For a prefix height row $h=(x_1,\ldots,x_d,z)$, the update is
\[
 hM_k=(x_2,\ldots,x_d,kx_1+z,(k-1)x_1+z).
\]
Starting from the row of ones, the invariant cone is
\begin{equation}\label{eq:cone}
 x_1\le\cdots\le x_d,\qquad z\le x_d\le x_1+z.
\end{equation}
Indeed, $x_d\le x_1+z\le kx_1+z$ gives the new ordering,
and the last inequality for the updated row reduces to
$x_1\le x_2$.

Put $B=\chi_2^{d-1}$. Direct iteration gives
\begin{align*}
 B(1)&=d,\\
 B(i)&=N\,1\,2\cdots(i-2)(i-1)^2\qquad(2\le i\le d),\\
 B(N)&=N\,1\,2\cdots(d-1).
\end{align*}
At an occurrence following $\theta=\theta_s$, write
\[
 U=\theta B(1)=\theta(d),\quad
 V_i=\theta B(i)\ (2\le i\le d),\quad W=\theta B(N),
 \qquad a=|U|=x_d.
\]
The cone gives the comparisons
\begin{equation}\label{eq:balance}
 \begin{split}
 a&\le |V_i|=z+\sum_{r=1}^{i-2}x_r+2x_{i-1}\le(i+1)a,\\
 a&\le |W|=z+\sum_{r=1}^{d-1}x_r\le da.
 \end{split}
\end{equation}
The lower bounds use $x_1+z\ge a$; the upper bounds use
$x_r,z\le a$.

After the full block $B\chi_KB$, the tower names are
\begin{align}
 Z_1&=WU^K,\notag\\
 Z_i&=WU^{K-1}V_2\cdots V_{i-1}V_i^2
       \qquad(2\le i\le d),\label{eq:Z}\\
 Z_N&=WU^{K-1}V_2\cdots V_d.\notag
\end{align}
Every name thus contains a power $U^{r_i}$ with $r_i\ge K-1$.
Using \eqref{eq:balance}, the length outside this power is at most
$C_da$, where
\begin{equation}\label{bruin:eq:outside-constant}
 C_d=d+\sum_{i=2}^{d-1}(i+1)+2(d+1)
     =\frac{d^2+7d-2}{2}=A_d-1.
\end{equation}
The displayed upper bound is largest for $Z_d$.

\subsubsection{Application of the BT tower estimates}

\begin{proof}[Proof of Proposition~\ref{prop:block}]
The parsing proof of Lemma~\ref{bt:lem:towers} is unchanged:
$2,\ldots,d$ are singleton images, and an $N$ followed by $k$
or $k-1$ letters $1$ identifies the remaining two images.
For $\theta=\mathrm{id}$, every word in \eqref{eq:Z} contains all
$N$ letters and has length at least $K+d$. Hence the selected
blocks have positive incidence and force the minimum height to
infinity, as in Lemma~\ref{bt:lem:growth}. In particular, the prefix
limits exist and $q_j\to\infty$.

In each $Z_i$-tower, discard the complement of $U^{r_i}$ and its
last copy of $U$. There are at most $A_da$ discarded levels and
at least $(K-2)a$ retained levels. The proof of
Proposition~\ref{bt:prop:mismatch} therefore gives
\begin{equation}\label{bruin:eq:one-symbol}
 \nu\{x:x_0\ne x_a\}\le\frac{A_d}{K+A_d-2}
\end{equation}
for every invariant symbolic probability. The same normalized return-word frequency argument used in
Proposition~\ref{bt:prop:coding}, with $C_d$ in place of $8$, gives
\begin{equation}\label{bruin:eq:frequencies}
 |\nu([w])-p_U(w)|\le\frac{C_da+m-1}{(K-1)a}
 \qquad(|w|=m).
\end{equation}
The compact symbolic system admits an invariant probability, and
\eqref{bruin:eq:frequencies} makes it unique. The correspondence in
\cite[Proposition~5 and Lemma~6]{Bruin} transfers this conclusion
and the discrepancy bound to the interval map. Lemmas~\ref{lem:generating}
and~\ref{bt:lem:rigidity} finish the proof exactly as in the BT case.
\end{proof}

\subsubsection{The induction cylinder for Bruin maps}
\label{bruin:sec:parameter}

The extension of the AHS measure in \cite{CS} has the same ergodic
full-support coding \cite[Sections~6.2 and 7.1]{AAFHS}.
The only change in the cylinder construction of
Section~\ref{sec:parameter} is that $d-1$ digits $2$ are required
on each side of $K$, as in \eqref{eq:block}.
The block \eqref{eq:block} can be realized inside an admissible
infinite Gauss itinerary by
\cite[equation~(2) and Lemma~4]{Bruin}.
As in \eqref{bt:eq:synchronized-path}, retain a finite marked
elementary word beginning at the switch that completes the
preceding Gauss digit and containing the entire block.
The cylinder defined by this word has positive measure by full support.
The recurrence argument and selection of disjoint occurrences with
$K_j\to\infty$ then apply verbatim. Proposition~\ref{prop:block}
completes the proof of Theorem~\ref{thm:typical}.

\section{A sufficient condition for nonrigidity}\label{sec:criterion}

We prove Theorem~\ref{thm:criterion} by comparing a candidate rigidity
time $q$ with a substitution level whose return heights are comparable
to $q$. Linear recurrence turns a small mismatch probability into
Hamming estimates for adjacent substituted words. Composing these
estimates along return words and using a signed relation between their
lengths produces an approximate period of fixed nonzero size.
Aperiodicity gives a uniform obstruction to this conclusion.
The identity for return spaces in \cite[Lemma~5.3]{Mercat}
allows us to use returns to one fixed letter, so that all rotations
in the signed combination act on the same prefix. We use
\cite[Proposition~5.2]{Mercat} separately in
Section~\ref{sec:bruin-weak-mixing} to exclude nontrivial measurable
eigenvalues for the Bruin construction.

\subsection{Recurrence and mismatch estimates}

Let $\sigma:\calA\to\calA^+$ be a primitive aperiodic substitution,
with language $\mathcal L_\sigma$, two-sided subshift $X_\sigma$,
shift $S$, and unique invariant probability $\nu$;
see \cite[Section~2.3.1]{DHS}.
A word is \emph{admissible} if it belongs to $\mathcal L_\sigma$.
For an admissible word $W$, write
$[W]=\{x\in X_\sigma:x_{[0,|W|)}=W\}$.
We use the column-count convention
\[
 M_{ij}=|\sigma(j)|_i,\qquad
 h_m=\one M^m,\qquad h_m(i)=|\sigma^m(i)|.
\]

\begin{lemma}
\label{nr:lem:PF}
There are constants $H>1$ and $K\in\N$ such that, for every
integer $q\ge1$, some $m=m(q)\ge0$ satisfies
\begin{equation}\label{eq:scale}
 q\le h_m(i)<Hq\qquad(i\in\calA).
\end{equation}
Moreover,
\begin{equation}\label{eq:recurrence}
 V,W\in\mathcal L_\sigma,\quad |W|\ge K|V|
 \quad\Longrightarrow\quad V\text{ occurs in }W.
\end{equation}
If $\ell M=\lambda\ell$ is a positive Perron--Frobenius row and
$C=(\max_i\ell_i)/(\min_i\ell_i)$, one may take
\[
 H=C^2\lambda,
 \qquad m(q)=\left\lceil\log_\lambda(Cq)\right\rceil.
\]
\end{lemma}

\begin{proof}
We have $\lambda>1$ since $\sigma$ is primitive and aperiodic.
Writing $\ell_-=\min_i\ell_i$ and $\ell_+=\max_i\ell_i$, the
coordinatewise inequalities
\[
 \frac{\ell}{\ell_+}\le\one\le\frac{\ell}{\ell_-}
\]
give, after multiplication by $M^m$,
\[
 C^{-1}\lambda^m\le h_m(i)\le C\lambda^m.
\]
The choice of $m(q)$ gives $Cq\le\lambda^{m(q)}<\lambda Cq$,
and hence \eqref{eq:scale}.
Primitive substitution subshifts are linearly recurrent
\cite[Proposition~25]{DHS}. Thus there is an integer $K_0$
such that consecutive occurrences of any admissible word of
length $L$ are at distance at most $K_0L$.
Taking $K=K_0+1$, every admissible word of length $KL$
contains every admissible word of length $L$, proving
\eqref{eq:recurrence}.
\end{proof}

For $q\ge1$ and an admissible word $W=W_0\cdots W_{L-1}$, put
\[
 \mathcal B_q=\{x\in X_\sigma:x_0\ne x_q\},\qquad
 d_\nu(q)=\nu(\mathcal B_q),
 \qquad H_q(W)=\#\{0\le j<L-q:W_j\ne W_{j+q}\}.
\]

\begin{lemma}\label{lem:local}
For every admissible word $W$ of length $L>q\ge1$,
\[
 H_q(W)\le(KL-q)d_\nu(q)\le KLd_\nu(q).
\]
\end{lemma}
\begin{proof}
By Lemma~\ref{nr:lem:PF}, the block $x_{[0,KL)}$ contains $W$
for every $x\in X_\sigma$. Write
$x_{[b,b+L)}=W$, where $0\le b\le KL-L$. Then
\[
 H_q(W)
 =\sum_{j=0}^{L-q-1}\mathbf1_{\mathcal B_q}(S^{b+j}x)
 \le\sum_{s=0}^{KL-q-1}\mathbf1_{\mathcal B_q}(S^s x),
\]
because the indices $b+j$ are distinct and lie between $0$ and
$KL-q-1$. Integrating and using shift invariance gives
\[
 H_q(W)
 \le\sum_{s=0}^{KL-q-1}\nu(S^{-s}\mathcal B_q)
 =(KL-q)d_\nu(q),
\]
as required.
\end{proof}

For a candidate time $q\ge1$, choose $m=m(q)$ as in
\eqref{eq:scale}, put $U_i=\sigma^m(i)$, and let $w_i$ be its
length-$q$ prefix. For words of length $q$, use the normalized
Hamming distance and cyclic rotations
\[
 \delta(v,w)=\frac1q\#\{0\le t<q:v(t)\ne w(t)\},\qquad
 (R^s w)(t)=w(t+s\bmod q).
\]
If $ij$ is an admissible pair, then $U_iU_j$ is admissible. The second block
begins at position $h_m(i)$, which determines the phase in the
following estimate.

The following lemma is the key estimate in our nonrigidity construction.
Its comparisons will be composed along return words and combined with
the algebraic height relation to obtain a fixed nonzero displacement.

\begin{lemma}\label{lem:edge}
If $ij$ is an admissible pair, then
\begin{equation}\label{eq:edge}
 \delta(w_j,R^{h_m(i)}w_i)\le2HK\,d_\nu(q).
\end{equation}
\end{lemma}
\begin{proof}
The word
$W=U_iU_j=\sigma^m(ij)$ is admissible and has length
$h_m(i)+h_m(j)<2Hq$. For each $0\le t<q$, write uniquely
\[
 h_m(i)+t=r_t+s_tq,
 \qquad 0\le r_t<q,\qquad s_t\ge1.
\]
The symbol of $W$ at position $h_m(i)+t$ is $w_j(t)$, while the
symbol at the residue position $r_t$ is
$w_i(r_t)=(R^{h_m(i)}w_i)(t)$.
Here the first $q$ positions of $W$ lie in $U_i$, because
$h_m(i)\ge q$.

Compare the symbols at successive positions in the arithmetic
progression
\[
 r_t,\ r_t+q,\ \ldots,\ r_t+s_tq=h_m(i)+t.
\]
If the symbols at its first and last positions differ, then
at least one consecutive pair has different symbols. Consequently
\[
 \mathbf1_{\{w_j(t)\ne(R^{h_m(i)}w_i)(t)\}}
 \le\sum_{\ell=0}^{s_t-1}
 \mathbf1_{\{W_{r_t+\ell q}\ne W_{r_t+(\ell+1)q}\}}.
\]
Every comparison on the right belongs to the index set defining
$H_q(W)$. Its lower endpoint is nonnegative, and its upper endpoint
is at most $h_m(i)+t\le h_m(i)+q-1<|W|$, where we used
$h_m(j)\ge q$. Thus no comparison leaves the admissible word.

The residues $r_t$ are pairwise distinct as $t$ ranges through
$0,\ldots,q-1$. Hence the positions $r_t+\ell q$, for
$0\le\ell<s_t$, are distinct: positions belonging to different
values of $t$ lie in different residue classes modulo $q$.
Summing over $t$ and applying Lemma~\ref{lem:local} gives
\begin{align*}
 q\,\delta(w_j,R^{h_m(i)}w_i)
 &\le H_q(W)\\
 &\le K|W|d_\nu(q)
 \le 2HKq\,d_\nu(q).
\end{align*}
\end{proof}

\subsection{Return words and a fixed displacement}

Composing the preceding comparisons along a return to the same
letter gives an approximate period of its prefix. A signed relation
between return lengths makes this displacement independent of the
substitution level.

\begin{proposition}\label{prop:signed}
Suppose that $r_1,\ldots,r_s$ are return words to one letter $a$,
that $n_1,\ldots,n_s\in\Z$, and that
\[
 u=\sum_{j=1}^s n_j\ab(r_j),\qquad
 Mu=\eta u,\quad\eta\in\{1,-1\},\qquad c=|\one u|>0.
\]
Put $B=\sum_{j=1}^s|n_j||r_j|$. Then
\begin{equation}\label{eq:quantitativegeneral}
 d_\nu(q)\ge\frac1{4BHK^2(c+1)}
 \qquad\bigl(q\ge2K(c+1)\bigr).
\end{equation}
\end{proposition}
\begin{proof}
For any $q\ge1$, take the corresponding prefixes $w_i$ and
height row $h_m=\one M^m$. Put
\[
 \varepsilon=2HKd_\nu(q),\qquad
 D_a(t)=\delta(w_a,R^t w_a).
\]
Fix $j$, write $r_j=a_0\cdots a_{\ell-1}$, and set
$a_\ell=a_0=a$ and $p_k=\sum_{t=0}^{k-1}h_m(a_t)$.
Since $r_ja$ is admissible, Lemma~\ref{lem:edge} applies to each
pair $a_ka_{k+1}$. Rotation invariance and the triangle inequality give
\[
 \delta(w_{a_{k+1}},R^{p_{k+1}}w_a)
 \le\varepsilon+\delta(w_{a_k},R^{p_k}w_a)
 \qquad(0\le k<\ell).
\]
Starting from $p_0=0$ and iterating, we return to the same prefix
$w_a$ with total phase
$p_\ell=h_m\ab(r_j)=|\sigma^m(r_j)|$. Hence
\begin{equation}\label{signed:eq:individual-return}
 D_a\bigl(h_m\ab(r_j)\bigr)\le |r_j|\varepsilon.
\end{equation}

Rotation invariance also gives $D_a(-t)=D_a(t)$ and
$D_a(t+t')\le D_a(t)+D_a(t')$. Thus the estimates for individual
returns can be combined with the signed coefficients $n_j$:
\[
 D_a(h_mu)
 =D_a\left(\sum_j n_jh_m\ab(r_j)\right)
 \le\sum_j|n_j|D_a\bigl(h_m\ab(r_j)\bigr)
 \le B\varepsilon.
\]
The eigenvector relation makes the absolute displacement independent
of the substitution level, since
\[
 h_mu=\one M^mu=\eta^m\one u\in\{c,-c\}.
\]
Thus $D_a(h_mu)=D_a(c)$, and therefore
\begin{equation}\label{eq:fixeddisplacement}
 D_a(c)\le2BHKd_\nu(q).
\end{equation}

By aperiodicity, choose an admissible word $V$ of length $c+1$ with
$V_0\ne V_c$, and put $L=K(c+1)$.
Each of the $\lfloor q/L\rfloor$ disjoint length-$L$ blocks
in $w_a$ contains $V$. They contribute distinct mismatches at
displacement $c$, with both endpoints inside $w_a$. Thus
\begin{equation}\label{signed:eq:fixed-lower-bound}
 D_a(c)\ge\frac{\lfloor q/L\rfloor}{q}\ge\frac1{2L}
 \qquad(q\ge2L).
\end{equation}
Combining the two bounds proves \eqref{eq:quantitativegeneral}.
\end{proof}

\begin{proof}[Proof of Theorem~\ref{thm:criterion}]
Put $\tau=\sigma^p$, so $M_\tau=M^p$ and $M_\tau u=u$.
Taking a power of a primitive substitution preserves its language
and subshift \cite[Section~2.3.2]{DHS}; hence
$(X_\tau,\nu,S)=(X_\sigma,\nu,S)$.
Fix a letter $a$, and let $\calR_a$ be the span of the count
vectors of its return words. Since $M_\tau u=u$, we have
$u=M_\tau^n u$ for every $n$, so $u$ belongs to the eventual
rational image
\[
 \mathcal I_\tau=\bigcap_{n\ge0}M_\tau^n\Q^{\calA}.
\]
By \cite[Lemma~5.3]{Mercat},
\[
 \mathcal I_\tau\cap\calR_\tau
 =\mathcal I_\tau\cap\calR_a.
\]
The equality of languages gives $\calR_\tau=\calR_\sigma$;
hence $u\in\calR_a$. This representation by returns to a single
letter is needed in Proposition~\ref{prop:signed}, where the signed
rotation estimates all concern the same prefix $w_a$.
Since $\calR_a$ is spanned by integer vectors and $u$ is rational,
clearing denominators gives an integer $D\ge1$, return words
$r_1,\ldots,r_s$ to $a$, and integers $n_1,\ldots,n_s$ such that
\[
 u_0=Du=\sum_{j=1}^s n_j\ab(r_j),\qquad M_\tau u_0=u_0.
\]
Set $c=|\one u_0|>0$ and $B=\sum_j|n_j||r_j|$.
Applying Proposition~\ref{prop:signed} to $\tau$ gives
\[
 d_\nu(q)\ge\frac1{4BH_\tau K_\tau^2(c+1)}
 \qquad\bigl(q\ge2K_\tau(c+1)\bigr),
\]
where $H_\tau,K_\tau$ are the constants of
Lemma~\ref{nr:lem:PF}. This proves \eqref{eq:gap}.
Finally,
\[
 \sum_{a\in\calA}
 \|\mathbf1_{[a]}\circ S^q-\mathbf1_{[a]}\|_{L^2(\nu)}^2
 =2d_\nu(q),
\]
so the positive lower bound excludes rigidity.
\end{proof}

\subsection{Weak mixing for the Bruin construction}
\label{sec:bruin-weak-mixing}

For self-similar Bruin maps, the algebraic condition also excludes
nontrivial measurable eigenvalues.

\begin{lemma}
\label{lem:bruin-nonrigidity-weak-mixing}
Let $T$ be a self-similar infinite-type Bruin interval translation
map whose return substitution $\sigma$ over one period of the Gauss
induction \eqref{eq:gauss} is primitive and aperiodic. If $\sigma$ satisfies the hypotheses of
Theorem~\ref{thm:criterion}, then $T$ is weakly mixing with respect
to its unique invariant probability.
\end{lemma}

\begin{proof}
Write $\calA=\{1,\ldots,N\}$, and let $M$ be the incidence matrix
of $\sigma$. Every elementary Bruin substitution in \eqref{eq:chi}
has a unimodular incidence matrix. Moreover, the first letter of
$\chi_k(a)$ is $\min(a+1,N)$, independently of $k$.
If the period consists of $\ell$ elementary substitutions, the first
letter of $\sigma^r(a)$ is $\min(a+\ell r,N)$.
Thus $M$ is unimodular, and, choosing $\ell r\ge N-1$, every image
under $\tau=\sigma^r$ begins with $N$.
For each $a$, choose $b$ such that $ab$ is admissible.
The word $\tau(a)N$ is a prefix of $\tau(ab)$, so
$\tau(a)$ is a concatenation of return words to $N$.
Hence the columns of $M^r$, which form an integer basis, belong
to the return lattice, proving
\begin{equation}\label{eq:bruin-full-return-lattice}
 \operatorname{span}_{\Z}
 \{\ab(w):w\text{ is a return word to }N\}=\Z^N.
\end{equation}
In particular, $\calC_\sigma=\calR_\sigma^\perp=0$.

Let $e^{2\pi it}$ be a measurable eigenvalue of the substitution
system. By \cite[Proposition~5.2]{Mercat} (see also \cite{Host}),
applied with the initial morphism equal to the identity, there are
$k\ge0$, a real coboundary $b\in\calC_\sigma$, and $w\in\Z^N$
such that
\[
 (t\one M^k+b-w)M^n\longrightarrow0.
\]
Since $\calC_\sigma=0$, we have $b=0$. Unimodularity then gives
$z=wM^{-k}\in\Z^N$, and hence
\begin{equation}\label{eq:bruin-stable-eigenvalue}
 v=t\one-z,\qquad vM^n\longrightarrow0.
\end{equation}
For the vector $u$ in Theorem~\ref{thm:criterion},
$vu=vM^{pn}u\to0$, so
\[
 t=\frac{zu}{\one u}\in\Q.
\]
Choose an integer $D\ge1$ with $Dv\in\Z^N$.
The integer rows $DvM^n$ tend to zero and hence are eventually
zero. Invertibility gives $v=0$, so $t\in\Z$.
The subshift associated with $\sigma$, and hence $T$, is weakly mixing.
\end{proof}

\section{The eight-interval example}\label{sec:eight}
\subsection{The map and its self-inducing interval}

Let $\varrho$ be the root of $Q$ specified in
Theorem~\ref{thm:eight}. Since $Q$ is strictly decreasing on
$(0,1/3)$ and $Q(1/4)>0>Q(1/3)$, we have
$1/4<\varrho<1/3$. Put $\kappa=2/(1-\varrho^7)$ and define
the interval lengths by
\begin{equation}\label{eq:geometric-lengths8}
 \ell=
 \left(\kappa\varrho^2,\kappa\varrho^4,\kappa\varrho^6,
       \kappa\varrho,\kappa\varrho^3,\kappa\varrho^5,
       \kappa\varrho^7,\frac{1-3\varrho}{1-\varrho}\right)^{\mathsf T}.
\end{equation}
They are positive, and
$\kappa\sum_{j=1}^7\varrho^j+(1-3\varrho)/(1-\varrho)=1$.
Set $b_0=0$, $b_i=\sum_{j=1}^i\ell_j$, and
$I_i=[b_{i-1},b_i)$. The tail sums $a_i=1-b_i$, $i\le7$,
agree with \eqref{eq:parameters}, using $Q(\varrho)=0$;
in particular, $a_4=\varrho$ and $a_7=\ell_8$.
Thus $T=T_a$ translates $I_i$ by $a_i$ for $i\le7$, and
$I_8$ by $a_7-1$.

We compute the first return to
\[
 J=[1-\varrho,1)=I_5\cup\cdots\cup I_8,
 \qquad \phi(x)=1-\varrho+\varrho x,
 \qquad J_i=\phi(I_i).
\]
The length formulas give the endpoint identities
\begin{equation}\label{eq:eight-return-domains}
 \begin{gathered}
 J_1=I_5,\quad J_2=I_6,\quad J_3=I_7,\\
 T(J_4)=I_1,\quad T(J_5)=I_2,\quad T(J_6)=I_3,\\
 T(J_7)=[b_3,1-3\varrho),\qquad
 T(J_8)=[1-3\varrho,a_7).
 \end{gathered}
\end{equation}
The two inequalities needed for the last branches are
\begin{equation}\label{eq:eight-endpoint-gaps}
 1-3\varrho-b_3=\varrho\ell_7>0,
 \qquad
 1-2\varrho-a_7=\frac{2\varrho^2}{1-\varrho}>0.
\end{equation}
On $I_4=[b_3,1-\varrho)$, the map adds $\varrho$.
Consequently the intervals
\[
 [b_3+j\varrho,1-(3-j)\varrho)\quad(0\le j\le2),
 \qquad
 [1-(3-j)\varrho,a_7+j\varrho)\quad(0\le j\le1)
\]
lie in $I_4$ and describe the successive visits of $J_7$ and
$J_8$ before return. Their next images lie in $J$, since
$b_3+3\varrho>1-\varrho$ and $a_7+2\varrho<1$.
For $J_4,J_5,J_6$, the first images are outside $J$ and the
second images are inside: $a_3>3\varrho>1-\varrho$ gives
$T(I_i)\subset J$ for $i=1,2,3$.
Finally, $\ell_{i+4}=\varrho\ell_i<\varrho$ for $i=1,2,3$,
so $T(J_i)\subset J$. The first-return words are therefore
\[
 (5,6,7,81,82,83,8444,844),
\]
exactly the substitution $\zeta$ in \eqref{eq:zeta}.

Its column-count incidence matrix is
\begin{equation}\label{eq:M8}
 M=\begin{pmatrix}
 0&0&0&1&0&0&0&0\\
 0&0&0&0&1&0&0&0\\
 0&0&0&0&0&1&0&0\\
 0&0&0&0&0&0&3&2\\
 1&0&0&0&0&0&0&0\\
 0&1&0&0&0&0&0&0\\
 0&0&1&0&0&0&0&0\\
 0&0&0&1&1&1&1&1
 \end{pmatrix}.
\end{equation}
The signed translation row $t=(a_1,\ldots,a_7,a_7-1)$ is a
left eigenvector of $M$ with eigenvalue $\varrho$:
\begin{equation}\label{eq:eight-translation-eigenvector}
 tM=\varrho t.
\end{equation}
The substitution acts here on translations; the interval lengths
are recovered as $\ell_i=a_{i-1}-a_i$, with $a_0=1$ and $a_8=0$.
The $i$th coordinate of $tM$ is the total translation along
$\zeta(i)$. Thus, for $x\in I_i$,
\[
 T_J(\phi(x))=\phi(x)+(tM)_i
             =\phi(x)+\varrho t_i=\phi(T(x)).
\]
This proves the self-induction identity
\begin{equation}\label{eq:eight-inducing}
 T_J\circ\phi=\phi\circ T.
\end{equation}

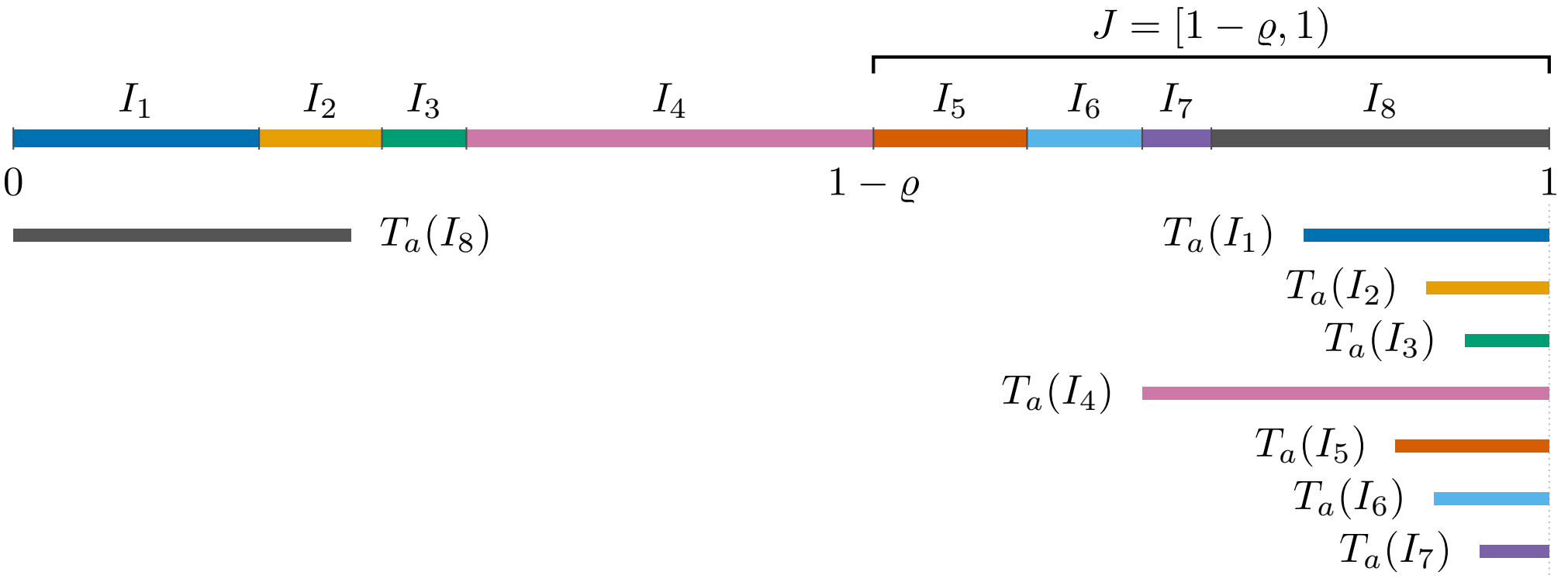
\begin{figure}[htbp]
\centering
\begin{tikzpicture}[x=13.1cm,y=1cm,font=\small,line cap=butt]
\definecolor{itmone}{HTML}{0072B2}
\definecolor{itmtwo}{HTML}{E69F00}
\definecolor{itmthree}{HTML}{009E73}
\definecolor{itmfour}{HTML}{CC79A7}
\definecolor{itmfive}{HTML}{D55E00}
\definecolor{itmsix}{HTML}{56B4E9}
\definecolor{itmseven}{HTML}{7B61A8}
\definecolor{itmeight}{HTML}{555555}
% Schematic widths preserve the order of the actual interval lengths.
% The same schematic width is used for each interval and its image.
\foreach \l/\r/\col/\i in {
  0/.16/itmone/1,.16/.24/itmtwo/2,.24/.295/itmthree/3,
  .295/.56/itmfour/4,.56/.66/itmfive/5,.66/.735/itmsix/6,
  .735/.78/itmseven/7,.78/1/itmeight/8}{
  \fill[\col] (\l,4.20) rectangle (\r,4.35);
  \draw[black!70,thin] (\l,4.18)--(\l,4.37);
  \node at ({(\l+\r)/2},4.60) {$I_{\i}$};
}
\draw[black!70,thin] (1,4.18)--(1,4.37);
\node[below] at (0,4.15) {$0$};
\node[below] at (.56,4.15) {$1-\varrho$};
\node[below] at (1,4.15) {$1$};
\draw[thick] (.56,4.83)--(.56,4.97)--(1,4.97)--(1,4.83);
\node at (.78,5.22) {$J=[1-\varrho,1)$};
\draw[black!25,densely dotted] (1,.55)--(1,3.73);
\foreach \l/\y/\col/\i in {
  .84/3.45/itmone/1,.92/3.00/itmtwo/2,
  .945/2.55/itmthree/3,.735/2.10/itmfour/4,
  .90/1.65/itmfive/5,.925/1.20/itmsix/6,
  .955/.75/itmseven/7}{
  \fill[\col] (\l,\y-.055) rectangle (1,\y+.055);
  \node[anchor=east,xshift=-3pt] at (\l,\y) {$T_a(I_{\i})$};
}
\fill[itmeight] (0,3.395) rectangle (.22,3.505);
\node[anchor=west,xshift=3pt] at (.22,3.45) {$T_a(I_8)$};
\end{tikzpicture}
\caption{The eight-interval Bruin map, drawn schematically.
Upper intervals and their lower images have matching colours;
overlapping images are separated vertically, as in
\cite[Figure~1]{AAFHS}. For $i\le7$, $T_a(I_i)=[1-\ell_i,1)$,
whereas $T_a(I_8)=[0,\ell_8)$. The shortest intervals are enlarged
for legibility; their exact lengths are \eqref{eq:geometric-lengths8}.
The bracket marks the self-inducing interval $J=I_5\cup\cdots\cup I_8$.}
\label{fig:eight-itm}
\end{figure}

The first parameters at the four Gauss steps are
$a_1,a_2/a_1,a_3/a_2,\varrho/a_3$.
Since $a_3>3\varrho>3/4$, the first three lie in $(1/2,1)$
and the last in $(1/4,1/3)$.
The new last numerators at the first three steps are
$\varrho^2,\varrho^2a_1,\varrho^2a_2$; they preserve strict
positivity and ordering. Equation~\eqref{eq:eight-translation-eigenvector}
gives
\[
 G_7^4(a)=\varrho^{-1}
 (a_5,a_6,a_7,a_1+a_7-1,a_2+a_7-1,a_3+a_7-1,
                     3\varrho+a_7-1)=a.
\]
Hence the Gauss digits are $(1,1,1,3)^\infty$.
The scale $\varrho$ is the reciprocal of the Perron--Frobenius
eigenvalue of the period matrix on interval lengths. Indeed, in the
notation of \cite[Lemma~24]{AAFHS}, the nonnegative matrix
$Z_8(1)^3 Z_8(3)$ carries the induced length vector $\varrho\ell$
back to $\ell$, and $\ell$ is strictly positive.
Figure~\ref{fig:eight-itm} shows the map and its inducing interval.

\subsection{Infinite type and the invariant system}
\label{sec:eight-system}

The strict periodic self-induction established above gives
$a\in\calE_7$ and hence infinite type, by the self-similarity
argument of Boshernitzan--Kornfeld \cite{BK}; for the present
family, see \cite[Proposition~5]{Bruin}.

Direct calculation gives
\begin{equation}\label{eq:char8}
 \det M=1,\qquad \det(sI-M)=(s+1)Q(s),\qquad M^9>0.
\end{equation}
By \cite[Proposition~5]{Bruin}, the attracting Cantor dynamics
is identified with the subshift associated with $\zeta$.
Since $M^9>0$, this subshift is minimal and uniquely ergodic
\cite{DHS}; see also \cite[Lemma~6]{Bruin}.
It is infinite and hence aperiodic.
Write $\mu$ for its nonatomic invariant probability. The continuity coding
identifies it with the substitution probability $\nu$, modulo
the null endpoint exceptions discussed in Section~\ref{sec:rigidity}.

\subsection{Return words and the height relation}\label{sec:eight-returns}

Let $h_m=(|\zeta^m(i)|)_{i=1}^8=\one M^m$.
The integer column
\begin{equation}\label{eq:u8}
 u=(1,1,1,-1,-1,-1,-1,2)^{\mathsf T}
 \qquad\text{satisfies}\qquad Mu=-u,\quad\one u=1.
\end{equation}
Thus the height relation needed in Section~\ref{sec:criterion} is
\begin{equation}\label{eq:height8}
 h_m(1)+h_m(2)+h_m(3)-h_m(4)-h_m(5)-h_m(6)-h_m(7)
              +2h_m(8)=(-1)^m.
\end{equation}
Take the return words to $8$
\[
 (r_1,\ldots,r_8)=(81,82,83,844,8444,8445,8446,8447),
 \qquad v_j=\ab(r_j).
\]
Each $r_j8$ occurs in $\zeta^8(8)$, and each $r_j$ contains
exactly one $8$. Their count vectors satisfy
\begin{align}
 u&=v_1+v_2+v_3+v_4+v_5-v_6-v_7-v_8,
       \label{eq:returncombination}\\
 B&=\sum_{j=1}^8|r_j|=25.\label{eq:B25}
\end{align}
The identity \eqref{eq:returncombination} explicitly gives
$u\in\calR_8$, the representation by returns to one letter used
in the application of \cite[Lemma~5.3]{Mercat} above.
Proposition~\ref{prop:signed} therefore applies with
$\eta=-1$, $c=1$ and $B=25$.

\subsection{The constants and the uniform gap}

The admissible adjacent pairs are
\begin{equation}\label{eq:E}
 E=\{18,28,38,44,45,46,47,48,58,68,78,81,82,83,84\}.
\end{equation}
Direct calculations give
\[
 \begin{gathered}
 \max_j(M^9)_{ij}\le6\min_j(M^9)_{ij}\quad(1\le i\le8),
 \qquad \max_i|\zeta^9(i)|=2023,\\
 \max_i h_m(i)\le6\min_i h_m(i)\quad(0\le m\le8),
 \end{gathered}
\]
and every $\zeta^9(i)$ contains every pair in $E$.
Multiplication by $M^9$ extends the height-ratio bound to all
$m\ge9$. Since each first image has length at most four,
the first level with $\min_i h_m(i)\ge L$ satisfies
\begin{equation}\label{eq:H24}
 L\le h_m(i)<24L\qquad(1\le i\le8),\qquad H=24.
\end{equation}
An admissible word of length $L$ lies in two consecutive level-$m$
blocks, hence in every level-$(m+9)$ block. Such a block has
length less than $2023\cdot24L=48552L$.
Every admissible word of length $97104L$ contains a complete block
at this level. We may therefore take
\begin{equation}\label{eq:Kexplicit}
 K=97104.
\end{equation}

The only equal-letter pair in $E$ is $44$. No first image
starts with $4$, and each has at most three consecutive $4$'s.
Thus every constant run has length at most three, and the
length-$q$ prefix $w_8$ satisfies
\[
 D_8(1)=\delta(w_8,Rw_8)
 \ge\frac{\lceil q/3\rceil-1}{q}\ge\frac16
 \qquad(q\ge6).
\]
Combining this with \eqref{eq:fixeddisplacement} gives
\[
 \frac16\le2\cdot25\cdot24\cdot97104\,d_\nu(q)
 \qquad(q\ge6).
\]
The continuity coding identifies $d_\nu(q)$ with
$\mu\{x:\calP(x)\ne\calP(T^qx)\}$.
Consequently \eqref{eq:explicitgap} holds with
\[
 \eps_8=\frac1{699\,148\,800},
\]
which completes the proof of Theorem~\ref{thm:eight}.

\section*{Acknowledgements}
We heartily thank ChatGPT (Astra model) for the computations that leaded
us to the construction of the example in Section 4 and certain improvements
to the presentation. The authors take full responsibility for the accuracy
and originality of the work.

AS gratefully acknowledges the hospitality of CIRM (Marseille) and the
financial support of Basic Research Program of HSE University
(HSE-BR-2025-062).

\bigskip
\noindent\begin{minipage}{\textwidth}
\small
\textbf{Pascal Hubert}\par\nobreak\smallskip
\noindent Aix-Marseille Universit\'e, Institut de Math\'ematiques de
Marseille (I2M), UMR 7373,\par\nobreak
\noindent 3 place Victor Hugo, Case 19, 13331 Marseille Cedex 3, France\par\nobreak
\noindent\textit{Email address:}
\href{mailto:pascal.hubert@univ-amu.fr}{\nolinkurl{pascal.hubert@univ-amu.fr}}

\medskip
\noindent\textbf{Alexandra Skripchenko}\par\nobreak\smallskip
\noindent Faculty of Mathematics, National Research University Higher School
of Economics, Usacheva St. 6, 119048 Moscow, Russia\par\nobreak
\noindent\textit{Email address:}
\href{mailto:sashaskrip@gmail.com}{\nolinkurl{sashaskrip@gmail.com}}
\end{minipage}

\end{document}